\documentclass[12pt,amsmath,amssym,latexsymb,leqno]{article}

\usepackage{amssymb,latexsym}

\usepackage{graphicx}

\usepackage{psfrag}

\usepackage[final,dvips]{epsfig}

\usepackage{verbatim}

\usepackage{amsfonts}

\def\newtext{}
\def\newwtext{}

\makeatletter
\def\@begintheorem#1#2{\trivlist \item[\hskip \labelsep{\sc #1\ #2.}]\it}
\def\@opargbegintheorem#1#2#3{\trivlist
      \item[\hskip \labelsep{\sc #1\ #2.\ (#3)}]\it}
\def\@sect#1#2#3#4#5#6[#7]#8{\ifnum #2>\c@secnumdepth
     \let\@svsec\@empty\else
     \refstepcounter{#1}\edef\@svsec{\csname the#1\endcsname.\hskip 1em}\fi
     \@tempskipa #5\relax
      \ifdim \@tempskipa>\z@
        \begingroup #6\relax
          \@hangfrom{\hskip #3\relax\@svsec}{\interlinepenalty \@M #8\par}%
        \endgroup
       \csname #1mark\endcsname{#7}\addcontentsline
         {toc}{#1}{\ifnum #2>\c@secnumdepth \else
                      \protect\numberline{\csname the#1\endcsname}\fi
                    #7}\else
        \def\@svsechd{#6\hskip #3\relax  
                   \@svsec #8\csname #1mark\endcsname
                      {#7}\addcontentsline
                           {toc}{#1}{\ifnum #2>\c@secnumdepth \else
                             \protect\numberline{\csname the#1\endcsname}\fi
                       #7}}\fi
     \@xsect{#5}}
\def\section{\@startsection {section}{1}{\z@}{-1.5ex plus-1ex minus
    -.2ex}{-2.5ex plus.2ex}{\reset@font\bf}}
\def\subsection{\@startsection{subsection}{2}{\z@}{-3.25ex plus-1ex
    minus-.2ex}{-1.5ex plus.2ex}{\reset@font\sl}}

\makeatother

\newtheorem{predf}{Definition}[section]

\newtheorem{prex}[subsection]{Example}

\newtheorem{theo}[subsection]{Theorem}
\newtheorem{definition}[subsection]{Definition}
\newtheorem{prop}[subsection]{Proposition}
\newtheorem{lemma}[subsection]{Lemma}

\begin{document}
\newcommand{\C}{{\Bbb C}}
\newcommand{\R}{{\Bbb R}}
\newcommand{\D}{{\Bbb D}}
\newcommand{\re}{{\rm Re \,}}
\newcommand{\im}{{\rm Im \,}}
\newcommand{\e}{\varepsilon}
\newcommand{\voltp}[1]{V^{(#1)}}
\newcommand{\currp}[1]{I^{(#1)}}
\newcommand{\forwp}[1]{U^{(#1)}}
\newcommand{\Rn}[1]{{\Bbb R}^{#1}}

\def \beq {\begin {eqnarray}}
\def \eeq {\end {eqnarray}}
\def \ba {\begin {eqnarray*}}
\def \ea {\end  {eqnarray*}}
\def \p {\partial}
\def \tilde {\widetilde}
\def \hat {\widehat}

     \newcommand{\be}{\begin{equation}}   
\newcommand{\ee}{\end{equation}}
\newcommand{\ben}{\[}                
\newcommand{\een}{\]}
\newcommand{\bea}{\begin{eqnarray}}  
\newcommand{\eea}{\end{eqnarray}}
\newcommand{\bean}{\begin{eqnarray*}}
\newcommand{\eean}{\end{eqnarray*}}
\parindent=.5cm

\newcommand{\F}{{\mathcal{F}}}
\newcommand{\eps}{\varepsilon}
\newcommand{\nbr}{\mathcal{N}}
\newcommand{\tr}{^{\rm T}}
\newcommand{\ith}{^{\rm th}}
\newcommand{\domain}{\Omega}
\newcommand{\compdomain}{\Omega_m}

\def\mtrx#1#2{
  \left(
    \begin{array}{#1}
      #2
    \end{array}
  \right)}

\title{Inverse conductivity problem with an imperfectly known boundary}

\author{
Ville Kolehmainen\footnote{University of Kuopio, Department of Applied Physics,
P.O.Box 1627, FIN-70211 Finland},$ $\quad Matti Lassas\footnote{
Helsinki University of Technology, Institute of Mathematics,
P.O.Box 1100, FIN-02015, Finland},$ $\quad  and 
Petri Ola\footnote{Rolf Nevanlinna Institute, University of Helsinki
Helsinki, P.O.Box 4, FIN-00014, Finland}
\\
}
\date{}
\maketitle

{\bf Abstract.} We show how to eliminate the error caused by an incorrectly 
modeled boundary in electrical impedance tomography (EIT). 
In practical measurements, one usually lacks the exact knowledge of 
the boundary. Because of this the numerical reconstruction from the
measured EIT--data is done using a model domain that represents
the best guess for the true domain. However, it has been noticed
that the inaccurate model of the boundary causes severe errors
for the reconstructions. We introduce a new algorithm to find a 
deformed image of the original isotropic conductivity based on the
theory of Teichm\"uller spaces and implement it numerically.


\bigskip

\section{Introduction.} 
We consider the electrical impedance tomography (EIT) problem, i.e.
the determination of an unknown conductivity distribution 
inside a domain, for example the 
human thorax, from voltage and current measurements made on the 
boundary. Mathematically this is 
formulated as follows: Let $\Omega\subset \R^2$ be the measurement domain, and 
denote by $\gamma = (\gamma ^{ij})$ the symmetric
matrix describing the conductivity in $\Omega$. 
We assume that the matrix has components in 
$L^\infty (\Omega)$ and that it is strictly positive definite,
that is, for some $c>0$ we have
 $\langle\xi,\gamma (x)\xi\rangle \geq c||\xi||^2$ for all $x\in \Omega$. 
The electrical potential $u$ satisfies in $\Omega$ the 
equation
\beq\label{joht}
\nabla\cdot \gamma\nabla u = 0.
\eeq
To uniquely fix the solution $u$ it is enough to give its value on the boundary. Let this be $u|_{\p \Omega}=f
\in H^{1/2}(\p \Omega)$ where $ H^{1/2}(\p \Omega)$
is the Sobolev space. Then (\ref{joht}) 
has a unique weak solution $u\in H^1(\Omega)$.

Our boundary data is the map that takes the voltage distribution $f$ on the boundary for all $f$ to the corresponding 
current flux through the boundary, $\nu\,\cdotp \gamma\nabla
u$, where $\nu$ is the exterior unit normal to 
$\Omega$. Mathematically this amounts to the knowledge of the
 Dirichlet--Neumann map $\Lambda$ 
corresponding to $\gamma$, i.e. the map taking the Dirichlet boundary values to the 
corresponding Neumann boundary values 
of the solution to (\ref{joht}),
\ba
\Lambda_\gamma:u|_{\p \Omega}\mapsto \sum _{i,j=1}^2
 \nu_i\, \gamma ^{ij}\frac{\partial u}{\partial x_j}\bigg |_{\p \Omega}.
\ea
This defines a bounded operator
 $\Lambda_\gamma:H^{1/2}(\partial \Omega)\to H^{-1/2}(\p \Omega)$.
The symmetric quadratic form  corresponding to $\Lambda_\gamma$,
\beq\label{joule}
\Lambda_\gamma[h,h]:=\int_{\p \Omega}h\,\Lambda_\gamma h\,dS=
\int_{\p \Omega}\nabla u\cdotp\gamma\nabla u \,dx
\eeq
equals in physical terms the power needed to maintain the potential
$h$ on $\p \Omega$.

When  $\gamma$ is a scalar valued function times identity matrix,
we say that the conductivity is {\em isotropic}. As usual, conductivities
that may be matrix-valued are referred as  {\em anisotropic} conductivities.
The EIT--problem is to reconstruct $\gamma$ from $\Lambda_\gamma$.
The problem was originally proposed by Calder\'on \cite{ca} and then solved in
 dimensions three and higher
for isotropic smooth conductivities in \cite{su}. The two dimensional case that 
is relevant to us was
solved by Nachman \cite{na2} for isotropic conductivities assuming $\gamma\in W^{2,p}$, $p>1$
and then finally for general $L^\infty$--smooth isotropic 
conductivities by Astala and P\"aiv\"arinta in a celebrated paper 
\cite{ap}.
\medskip

The conductivity equation is invariant under deformations of the domain 
$\Omega$ in the following sense. If $F$ is a 
diffeomorphism taking $\Omega$ to some other domain
 $\tilde \Omega $, then $u\circ F^{-1}$ will satisfy the conductivity
equation in $\tilde \Omega $ with conductivity 
\beq\label{pushforward}
\tilde \gamma (x) = \left.\frac{F'(y) \, \gamma(y)\, (F'(y))^t}{|\det F'(y)|}
\right|_{y= F^{-1}(x)},
\eeq
where $F'$ is the Jacobi matrix of map $F$
and $u$ is a solution of $\nabla\cdotp \gamma\nabla u=0$ in $\Omega$.
We say that $\tilde \gamma$ is the push forward of $\gamma$ by $F$
and denote it by $\tilde \gamma=F_*\gamma$.
Note that all this is well defined for general matrix
valued $\gamma$. For us the starting point is the trivial observation
 that even if $\gamma$ is isotropic, the deformed
 conductivity $\tilde \gamma $ will not in general be isotropic.
The boundary measurements are invariant: When
$f:\p \Omega\to \p\tilde \Omega$ is the restriction of $F:\Omega\to 
\tilde \Omega$,
we say that $\tilde \Lambda = f_*\Lambda_\gamma$,
\[
((f_*\Lambda_\gamma)h)(x)= \left.(\Lambda_\gamma(h\circ f))(y)
\right|_{y= f^{-1}(x)},\quad h\in H^{1/2}(\p \tilde\Omega)
\]
is the push forward of $\Lambda_\gamma$ in $f$. 
As seen in \cite{sy}, it turns out that $f_*\Lambda_ \gamma=\Lambda_{F_*
\gamma}
$. 

The facts that the anisotropic conductivity equation and the boundary
measurements are invariant has the
 important consequence that the EIT--problem with an anisotropic 
conductivity is not uniquely solvable, even though the
isotropic problem is, see \cite{sy}.
\medskip

In practice when solving the EIT--problem in a given domain $\Omega$,
one typically 
seeks for the isotropic conductivity that minimizes
\beq\label{minimization problem 1}
||\Lambda_{meas} - \Lambda _{\gamma}||^2 +\alpha ||\gamma||_X ^2
\eeq
for $\gamma$ defined in terms of some triangulation of $\Omega$ 
as e.g. a piecewise constant function and $||\cdot||_X$ is some 
regularization norm 
\cite{kaipio00a}.
Here $\Lambda_{meas}$ is the measurement of the Dirichlet--Neumann map 
that contains measurement errors. 

In practice, one of the key difficulties in solving the EIT problem is
that the domain $\Omega$ may not be known accurately. It has been 
noticed that the use of slightly incorrect
model for $\Omega$, i.e., a slightly incorrect model of the boundary 
causes serious errors in reconstructions, see e.g. 
\cite{kolehmainen97e,adler96b,gersing96}. 
As an example, consider the EIT measurements of pulmonary function 
from the human thorax. The measurement electrodes are attached on 
the skin of the patient around the thorax. In principle, an exact
parameterization for the shape of the thorax could be obtained from other
medical imaging modalities such as magnetic resonance imaging (MRI) or 
computerized tomography (CT). However, in most cases 
these data is not available, and one has to resort to some approximate thorax
model. Further, the shape of the thorax varies between
breathing states, and it is also dependent on the orientation of 
the patient. Thus, the thorax geometry is known inaccurately 
even in the best case scenarios.

In this paper our 
aim is to propose a method to overcome the problem that the boundary
and its parameterization are not exactly known.
The set--up of the problem we consider is the following: 

We want to recover 
the unknown conductivity $\gamma$
in $\Omega$ from the measurements of Dirichlet-to-Neumann map, 
and we  assume a'priori that $\gamma$ is isotropic.
We assume $\partial \Omega$ and $\Lambda_\gamma$ are not known. 
Instead, let $\Omega_m$, called the model domain,
be our best guess for the domain and let
 $f_m:\partial \Omega \to\partial \Omega_{m}$ 
be a diffeomorphism modeling the approximate knowledge of the boundary.
As the data for the inverse problem, we assume that we 
are given the boundary of the model domain,
$\partial \Omega_m$ and the map  
$ \Lambda_m:=(f_m)_*\Lambda_\gamma$ on $\partial \Omega_m$.
Note that we have simplified the problem by assuming
that the only error in $\Lambda_m$ comes from the imperfect
knowledge of the boundary.

This set--up is motivated by the fact that the quadratic form
corresponding to $\Lambda_m$,
\ba
\Lambda_m[h,h]=\int_{\p \Omega_m}h\,\Lambda_m h\,dS=
\int_{\p \Omega}(h\circ f_m)\,\Lambda_\gamma (h\circ f_m)\,dS
,\quad\quad
h\in H^{1/2}(\p \Omega_m)
\ea
 equals physically to the power needed to maintain the potential $h\circ f_m$ on $\p \Omega$. 

Since $\Lambda_m$ usually does not correspond to
any isotropic conductivity because of the deformation done when going from the 
original domain $\Omega$ to $\Omega_m$, we obtain an erroneous solution
$\gamma$ when solving the minimization problem
(\ref{minimization problem 1}).
This means that a systematic error in domain model
causes a systematic error to the reconstruction.
In particular, local changes  
of the conductivity often give raise to non-localized changes in reconstructions
due to the above  systematic error.
Thus the spatial resolution of details of reconstructions are often weak. 
This is clearly seen in practical  measurements, see e.g. \cite{gersing96}.

We note that one could forget 
in  solving of the minimization problem
(\ref{minimization problem 1}) the assumption that $\gamma$ is
isotropic, and find the minimizer in the set of anisotropic conductivities.
However, the anisotropic inverse problem has non-unique solution,
and as the minimization problem is highly non-convex, the minimization
would be hard and, as usual, forgetting existing a'priori information makes
the solution significantly worse.

To formulate our main results, let us define certain  concepts. 
We start with the  maximal anisotropy of an anisotropic
conductivity.

\begin{definition}
Let $ \gamma^{jk}(x)$ be an $L^\infty(\Omega)$-smooth matrix valued conductivity in  $\Omega$
and let $\lambda_1(x)$ and  $\lambda_2(x)$, $\lambda_1(x)\geq
  \lambda_2(x)$ be the   
eigenvalues of matrix  $ \gamma^{jk}(x)$. 
We define  the maximal anisotropy of a conductivity to be
$K(\gamma)$ given by 
\ba
K(\gamma)= \sup_{x\in \Omega}K(\gamma,x),\quad\hbox{where}\ \ K(\gamma,x)=
\frac {\sqrt {L(x)}-1}{\sqrt {L(x)}+1},\quad {L(x)}=
\frac{\lambda_1(x)}{\lambda_2(x)}.
\ea
We call the function $K(\gamma,x)$ the anisotropy of $\gamma$ at $x$.
{\newtext Here $\sup$ denotes the essential supremum.}
\end{definition}

Sometimes, to indicate the domain $\Omega$,
 we denote $K(\gamma)=K_\Omega(\gamma)$.
As a particularly important example needed later, let us
 consider the conductivity matrices of the form
\beq\label{MAIN}
\hat \gamma(x)=
\eta(x) R_{\theta(x)}
\left(\begin{array}{cc}
\lambda^{1/2} & 0\\
0 & \lambda^{-1/2} \end{array}\right)  R_{\theta(x)}^{-1}
\eeq
where $\lambda\geq 1$ is a constant, 
$\eta(x)\in \R_+$ is a real valued function,
$R_{\theta(x)}$ is a rotation matrix corresponding 
to angle $\theta(x)$, where
\ba
 R_{\theta}=
\left(\begin{array}{cc}
\cos \theta  & \sin \theta\\
-\sin \theta & \cos \theta \end{array}\right).
\ea 
We denote such conductivities by $\hat \gamma=\hat \gamma_{\lambda,\theta,
\eta}$.
These conductivities have the anisotropy $K(\hat \gamma,x)=c_\lambda=(\lambda^{1/2}-1)/(\lambda^{1/2}+1)$ at every point and
thus their maximal anisotropy is $K=c_\lambda$. We call
such conductivities $\hat \gamma$ {\em uniformly anisotropic conductivities}.

The following theorem is the main result of the paper.

\begin{theo}\label{seainoaoikea}
Let $\Omega$ be a bounded, simply connected $C^{1,\alpha}$--domain with 
$\alpha >0$. Assume 
that $\gamma \in C^{0,\alpha}(\overline
\Omega)$ is an isotropic conductivity and $\Lambda_\gamma$ its Dirichlet--Neumann map. 
Let
$\Omega _m$ be a model of the domain (which is assumed to satisfy the same regularity 
assumptions as $\Omega$), and 
$f_m :\p\Omega\to \p\Omega _m$ be  a $C^{1,\alpha}$--smooth diffeomorphism.

Assume that we are given $\p \Omega_m$ and $\Lambda_m=(f_m)_*\Lambda_\gamma$.
Then 
\begin{enumerate}

\item There are unique 
$\lambda\geq 1$, $\theta\in L^\infty(\Omega_m,S^1)$ and $\eta
\in L^\infty(\Omega_m,\R_+)$
 such that the 
conductivity $\hat \gamma=\hat \gamma_{\lambda,\theta,\eta}$
satisfies $\Lambda_{\hat \gamma}=\Lambda_m$.

\item If $\gamma_2$ is an anisotropic conductivity in $\Omega_m$
such that  $\Lambda_{\gamma_2}=\Lambda_m$ then
$K(\gamma_2)\geq K(\hat \gamma)$. Moreover, 
$K(\gamma_2)= K(\hat \gamma)$ if and only if $\gamma_2=\hat \gamma$.
\end{enumerate}
\end{theo}

 Theorem \ref{seainoaoikea} can be interpreted by saying
that we can find a unique conductivity in $\Omega_m$ that is 
as close as possible to being  isotropic.

The proof of Theorem \ref{seainoaoikea} is based on the theory
of quasiconformal maps. There are several 
equivalent definitions for these maps, and we will present the one based on 
a partial differential equation (Beltrami equation) in Section 2.
However, the quasiconformal maps have 
also a geometric definition. Indeed, they are generalizations
of  conformal maps
that take infinitesimal disks at $z$ to
infinitesimal disks at $f(z)$, and the radii gets dilated 
by $|f'(z)|$. Analogously, a homeomorphic map is quasiconformal on a domain
$\Omega$ if infinitesimal disks at any $z\in \Omega$ get mapped to 
infinitesimal ellipsoids at $f(z)$. 
The ratio of the larger
semiaxis to smaller semiaxis is called the dilation of $f$ at $z$,
and the supremum of dilatations over $\Omega$
is the maximal dilation. 
This dilatation of infinitesimal discs is in fact the
reason why isotropic conductivities change to anisotropic ones
in push forwards with quasiconformal maps.

The crucial fact that we use proving 
Theorem \ref{seainoaoikea} is a result of Strebel \cite{st}, 
that roughly speaking says that  among all quasiconformal self--maps of the unit disk
to itself with a given sufficiently smooth boundary value there is a 
unique one with the minimal maximal dilation. This will yield 
that corresponding to the given 
boundary modeling map $f_m:\p \Omega\to \p \Omega_m$
there is a unique map $F:\Omega\to  \Omega_m$ having the minimal
maximal dilation. 
We will show that this leads to the following result:

\begin{prop}\label{ainoatoikeat}
Let $\Omega,$ $\gamma,$ $\Omega _m$, and $f_m$ satisfy
assumptions of Theorem \ref{seainoaoikea}.
Then there is  a unique map $F:\Omega\to \Omega_m$, depending
only on $f:\p \Omega\to \p \Omega_m$ such that
for the  uniformly anisotropic conductivity 
$\hat \gamma$ corresponding to $\gamma$ in 
Theorem \ref{seainoaoikea}
we have
\beq\label{determinantti}
\det (\hat \gamma(x))^{1/2}=\gamma(F^{-1}(x)).
\eeq
\end{prop}
Proposition \ref{ainoatoikeat} can interpreted as saying
that solving the minimization problem (\ref{minimization problem 1})
with $\alpha=0$
in the class of conductivities $\hat \gamma_{\lambda,\theta,\eta}$
we can find the function $(\det \hat \gamma(x))^{1/2}$ in 
$\Omega_m$ that represents a deformed image of original
conductivity $\gamma$ in the unknown domain $\Omega$
and the deformation depends only on the error made in modeling
the boundary, not on the conductivity in $\Omega$.

In particular, this turns out to be useful as local
perturbations of conductivity remain local in reconstruction: if 
we consider one fixed boundary modeling map $f_m:\p\Omega
\to \p\Omega_m$ but two
isotropic conductivities $\gamma_1$ and $\gamma_2=\gamma_1+\sigma$ 
in $\Omega$, then the reconstructions 
$\hat \gamma_1$ and $\hat \gamma_2$
obtained by Theorem \ref{seainoaoikea} corresponding
to $\gamma_1$ and $\gamma_2$ satisfy
\ba
\det(\hat \gamma_2(x))^{1/2}-\det(\hat \gamma_1(x))^{1/2}=\sigma(F^{-1}(x)).
\ea

{\bf Remark 1.} Theorem \ref{seainoaoikea} holds for
anisotropic conductivities in the sense that for
each $C^{1,\alpha}(\overline \Omega)$-smooth anisotropic
 conductivity $\gamma$ there is a unique
conductivity $\hat \gamma=\hat \gamma_{\lambda,\theta,\eta}$ such
that $(f_m)_*\Lambda_\gamma=\Lambda_{\hat \gamma}$. However,
for anisotropic conductivities Proposition \ref{ainoatoikeat}
does not hold as the map $F$ depends on the conformal class
of $\gamma$.
\medskip

The paper is organized as follows.
In Section 2 we consider isotropization of anisotropic 
conductivity using a diffeomorphism and pay close attention
to the smoothness required from $\gamma$ and $\Omega$ and 
introduce the necessary background from the theory 
anisotropic inverse problems.
We apply this in 
Section 3 in proving main results using the 
existence of a Teichm\"uller mapping. 
In Section 4 we consider physically realistic measurements,
i.e., so-called complete electrode model.
The numerical implementation for the complete electrode model 
is then described in the last sections. 
\bigskip

\section{Quasiconformal maps
and solvability of inverse problem with
anisotropic conductivity.} 
\label{sect: 2}
It is a classical
result that every Riemannian surface is locally conformal to
 a Euclidean plane: this corresponds to choosing the coordinate system to be {\em isothermal}. Similarly, every anisotropic conductivity matrix
can be transformed to an isotropic conductivity.
 We identify the
plane $\R^2$ with complex plane $\C$.

\begin{lemma}\label{isotropisointi}
Let $\Omega$ be a bounded, simply connected $C^{1,\alpha}$--domain with 
$\alpha >0$. Assume 
that $\gamma \in C^{0,\alpha}(\overline
\Omega)$ is an isotropic conductivity.
Then there is a $C^{1,\alpha}$-smooth
diffeomorphism $F:\overline \Omega\to \overline{\tilde \Omega},\ 
\tilde \Omega
=F(\Omega)\subset \C$
such that
\begin{equation}\label{epälinbeltrami}
F_*\gamma = \beta,
\end{equation}
where $F_*\gamma$ is defined by (\ref{pushforward}), and $\beta$
 is the identity matrix
multiplied by a $C^{0,\alpha}$-smooth scalar function. Moreover,
\[
\beta = (\det \gamma \circ F^{-1})^{1/2}I.
\]
\end{lemma}

The proof of this result is well known, but as
smoothness of $F$ is crucial later, we give the proof
for the convenience of the reader.

{\bf Proof.}
The equation (\ref{epälinbeltrami})
a'priori a nonlinear system 
for the derivatives of $F$.
However, in two dimensions this equation completely linearizes,
 and is equivalent to
the {\em Beltrami--equation}
\begin{equation}\label{beltrami}
\overline\partial F = \mu \partial F,
\end{equation}
where the complex derivatives are $\partial = \frac 12(
\frac \partial {\partial x}  - i\frac \partial {\partial y})$, 
$\overline\partial = \frac 12 (\frac \partial {\partial x}  + 
i\frac \partial {\partial y})$ and 
{\em the Beltrami coefficient} $\mu=\mu_F(z)$, called also
{\em complex dilatation} is given by
\beq\label{mu equation}
\mu = \frac {-\gamma_{11}+\gamma_{22}-2i\gamma_{12}}
{\gamma_{11}+\gamma_{22}+2\sqrt{\gamma_{11}\gamma_{22}-\gamma_{12}^2}}.
\eeq
The function $\mu$ has the crucial property that it is strictly less 
than one in modulus:
\begin{equation}\label{bound}
\sup_{z\in\Omega}|\mu (z)|<1.
\end{equation}
Let us extend the conductivity matrix $\gamma$
 (a'priori only defined in $\Omega$)
to the whole plane to be the identity matrix outside $\Omega$.
Similarly, $\mu$ is extended outside $\Omega$ by zero.

Next we consider how to solve the Beltrami equation, 
and for this we consider it in the whole
plane. In order to have a unique solution we fix the behaviour of $F$ at
 infinity.
Thus, consider
\beq \label{beltram}
& &\overline \p F(z)=\mu(z)\p F(z)\quad \hbox{in }\C, \\ \nonumber
& &F(z)=z+h(z),\\ \nonumber
& & \lim_{z\to\infty}h(z)=0
\eeq
where $\mu$ is a compactly supported $L^\infty$--function satisfying (\ref{bound}).
 This problem
has unique solution $F\in L^p_\delta$ when $p$
 is close enough to
2 and $-2/p<\delta<1-1/p$. For the proof of 
this see, for example \cite{ahl} or \cite{sy}.
 The proof is based on the fact that (\ref{beltram}) 
can be written
as an integral equation
\beq \label{beltram2}
F(z)-\frac 1{2\pi i}\int_\C \frac {\mu(\zeta)\p 
F(\zeta)}
{z-\zeta}\,da(\zeta)=z
\eeq
where $da(\zeta)$ is Euclidean area in $\C$ (or $\R^2)$. 
As $||\mu||_\infty <1$,  it turns out that
the left hand side of equation  (\ref{beltram2}) is of
the form of the identity plus a contractive operator in 
Sobolev space $W^{1,p}(\Omega)$, {\newtext with appropriate $p$,} and thus equation  (\ref{beltram2})
can be solved by an application of the Neumann--series argument.

Using interior Schauder estimates for equation 
(\ref{beltram}), we see that if $\gamma$ and thus $\mu$
are $C^{0,\alpha}$--smooth, the solution $F$ has 
to be locally  $C^{1,\alpha}$--smooth in $\C$, in particular
in $\overline \Omega$.
Using formula (\ref{pushforward}) we see that
$F_*\gamma$ is   $C^{0,\alpha}$--smooth
in closure of ${\tilde \Omega}$.
\hfill  $\Box$.

In general, any solution $F:\Omega\to \tilde \Omega$
to Beltrami equation for $\mu$ satisfying (\ref{bound})
and for which $F\in H^1(\Omega)$ 
is called {\em quasiregular}. If a quasiregular map 
 $F:\Omega\to \tilde \Omega$ is homeomorphism, 
it is said to be  {\em quasiconformal}.
The quasiconformality can be defined also in geometrical 
terms, see \cite{ahl}, \cite{le}.

Next we recall the recent results for inverse problems for
anisotropic conductivities $\gamma$.
 Let us consider a 
class of conductivities in $\Omega$, given by
\ba
\Sigma(\gamma)=\{ F_*\gamma\ | \ F:\Omega\to \Omega \hbox{ is homeomorphism
, }F,F^{-1}\in H^1(\Omega;\C), F|_{\p \Omega}=I\},
\ea
that is, $\Sigma(\gamma)$ is 
the equivalence class of the conductivity $\gamma$ in push forwards with
boundary preserving diffeomorphisms. Then 
$\Lambda_{\sigma}=\Lambda_{\gamma}$ for all
$\sigma\in \Sigma(\gamma)$. By \cite{alp}, the
converse is true, that is, if $\sigma$ is a
strictly positive definite $L^\infty$-conductivity
and $\Lambda_{\sigma}=\Lambda_{\gamma}$,
then $\sigma\in \Sigma(\gamma)$. In other
words, $\Lambda_{\gamma}$ determines the equivalence class
$\Sigma(\gamma)$.
Note that diffeomorphism 
 $F:\Omega\to \Omega$ such that $F\in H^1(\Omega;\C)$ 
and $F|_{\p \Omega}=I$ is quasiconformal.

\medskip

\section{Proof of main results}

We start with the proof of Theorem \ref{seainoaoikea}.
Note the fact that the conductivity $\gamma$ in $\Omega$ is isotropic
will not be used in the proof at all.
First we show that we can assume that $\Omega_m$
is the unit disc $\D\subset \C$. 

To prove this, let $f_m:\p \Omega\to \p \Omega_m$ be
the boundary modeling map and $\gamma$
a $C^{1,\alpha}$--smooth (also possibly anisotropic) conductivity in $\overline \Omega$.

Our first observation is that as 
$\Omega_m$  is a simply connected domain, it follows from
Riemann mapping theorem that it can be mapped to unit
disc $\D$ conformally. Moreover, as $\Omega_m$ is $C^{1,\alpha}$--smooth domain
it follows from Kellog-Warschawski theorem  \cite[Thm. 3.6] {pom},
that the Riemann--map can be chosen to be 
be 
a $C^{1,\alpha}$--diffeomorphism $F_0:\overline \Omega_m\to \overline \D$
such that $F_0:\Omega_m\to \D$ is conformal.
Thus, if $\sigma$ is some $C^{0,\alpha}$--smooth anisotropic
conductivity in $\overline\Omega_m$, we have that $\sigma_0=(F_0)_*\sigma$
is $C^{0,\alpha}(\overline \D)$--smooth conductivity.

Second, we observe that 
the  uniformly anisotropic conductivity $\hat \gamma_{\lambda,\theta,\eta}$ 
of the form (\ref{MAIN})
in $\Omega_m$ changes under $(F_0)_*$ to 
a uniformly anisotropic conductivity 
 $(F_0)_*\hat \gamma_{\lambda,\theta,\eta}
=\hat \gamma_{\lambda,\theta_0,\eta_0}$ in $\D$
such that 
$\eta_0= \eta\circ F^{-1}_0$.

Third, we see that as $F_0:\Omega_m\to \D$ is conformal,
the maximal anisotropy of $(F_0)_*\sigma$ and $\sigma$ satisfy
\ba
K_\D((F_0)_*\sigma)=K_\Omega(\sigma),
\ea
that is, the maximal anisotropy is preserved in conformal
transformations for any $\sigma$.

{\newtext Fourth, if $f_0=F_0|_{\p \Omega_m}$ then
 $\Lambda_{\sigma_0}=(f_0)_*\Lambda_\sigma$.
Also, we 
see that the our data is invariant in the change of the
model domain in the sense that 
$(\tilde f_m)_*\Lambda_\gamma=(\tilde f_0)_*((\tilde f_m)_*\Lambda_\gamma)$
where $\tilde f_m=f_0\circ f_m:\p \Omega\to \p \D$.}

These four observations yield that it is enough to prove
the assertion in the case when $\Omega_m=\D$. Indeed,
changing $\Omega_m$ to $\D$ with $F_0$ keeps 
the boundary measurements, the smoothness of objects, 
the maximal anisotropy as
well as class of uniformly anisotropic conductivities invariant.
More precisely, we can replace the boundary modeling 
map $f_m$ by the map $\tilde f_m=f_0\circ f_m$.

Thus, let us return proving Theorem \ref{seainoaoikea} in the 
case when $\Omega_m=\D$.
Let $f_m:\p\Omega\to \p \D$ be the boundary modeling
map that is a $C^{1,\alpha}$--smooth diffeomorphism
and $\gamma \in C^{0,\alpha}(\overline
\Omega)$ be an isotropic conductivity with Dirichlet--Neumann map
$\Lambda_\gamma$.

\begin{figure}[htbp]
\begin{center}
\psfrag{1}{$\Omega,\gamma$}
\psfrag{a}{$F_m$}
\psfrag{b}{$f_m$}
\psfrag{c}{$\D,\gamma_0$}
\psfrag{2}{$\Omega_1,\gamma_1$}
\psfrag{3}{$\D,\gamma_3$}
\psfrag{4}{$\D,\sigma_e$}
\psfrag{5}{$F_1$}
\psfrag{6}{$F_2$}
\psfrag{7}{$F_3$}
\psfrag{8}{$F_e$}
\psfrag{9}{$f_4$}

\includegraphics[width=12cm]{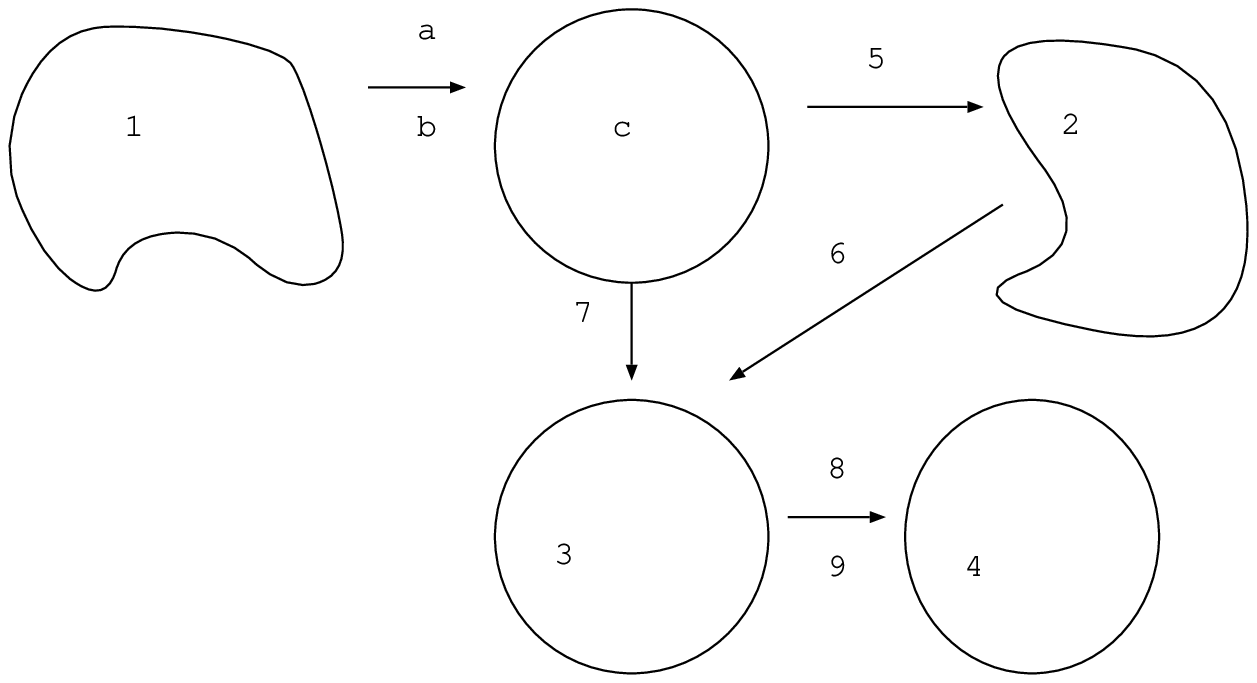}
\end{center}
\end{figure}

Let $F_m$ be some $C^{1,\alpha}(\overline \Omega)$--diffeomorphism
$F_m:\overline \Omega\to \overline \D$ such that
$F_m|_{\p \Omega}=f_m$. 
{\newwtext There are many ways to  construct such map, and 
for convenience of the reader, we present one simple way.
Let  $G: \Omega\to \D$ be a Riemann-map.
By \cite[Thm. 3.7]{pom}, $G$ has $C^{1,\alpha}$-extension 
$G:\overline \Omega\to \overline \D$. Let
$\phi= f_m\circ G^{-1}:\p \D\to \p \D$ and  
$\Phi(z)=|z|\exp(i(|z|^3\hbox{arg}(\phi(z/|z|))+
(1-|z|^3)\hbox{arg}(z/|z|)))$ be 
a $C^{1,\alpha}$-diffeomorphism $\D\to \D$ satisfying $\Phi|_{\p \D}=\phi$.
Then $F_m$ can be chosen to be the map $\Phi\circ G$.}

Let $\gamma_0=(F_m)_*\gamma$
be an anisotropic conductivity in $\D$. 
By Lemma \ref{isotropisointi}, there is a $C^{1,\alpha}$-diffeomorphism 
$F_1:\overline \Omega\to
\overline \Omega_1$ such that the conductivity $\gamma_1=(F_1)_*\gamma_0$
is $C^{0,\alpha}$--smooth isotropic conductivity
in the closure of the $C^{1,\alpha}$--smooth domain $\Omega_1$.

As $\Omega_1$  is a simply connected $C^{1,\alpha}$--smooth 
domain, by
 Kellog-Warschawski theorem cited above there is a
conformal map
$F_2: \Omega_1\to \D$ 
such that $F_2:\overline \Omega_1\to\overline \D$
is a $C^{1,\alpha}$--diffeomorphism. 
Let $F_3=F_2\circ F_1:\D\to \D$ and
$f_3=F_3|_{\p D}$.
Note that $(F_3)_*\gamma_0$ is isotropic conductivity in 
$\D$ as $F_2$ is conformal.

The boundary values of quasiconformal maps $\D\to \D$ are 
characterized as being the quasi-symmetric maps, that is, 
{\newtext homeomorphic maps  $f:\p\D\to\p \D$ such that 
$\theta (u) =\arg f(e^{iu})$ satisfy}
\begin{equation}\label{q-s}
k^{-1}\leq \frac{\theta(v)-\theta (u)}{\theta(v)-\theta(u)}
\leq k,\quad \hbox{for all }\ u,v\in \p \D,
\end{equation}
{\newtext with some $k>0$,}
see \cite{le}.

Let us consider next the map $f_4=f_3^{-1}:\p \D
\to \p \D$. 
Since $f_3$ and $f^{-1}_3$ are   $C^{1,\alpha}$--smooth, we see that
$f_4$ satisfies
\begin{equation}\label{strebe}
\lim _{t\to 0}\frac{\theta(u+t)-\theta (u)}{\theta(u-t)-\theta(u)}
 = 1\ {\rm uniformly\ in}\ u\in \p \D,
\end{equation}
and is in particular  quasi-symmetric. Thus $f_4$ is boundary value 
of at least one quasiconformal map. 
What is more, since $f_4$ satisfies
condition (\ref{strebe}), {\newtext it follows 
from results  Strebel \cite{st},
that
among all quasiconformal maps having $f_4$ as a boundary value
 there is a unique {\em extremal}
map $F_e$ in the sense that the $L^\infty$--norm of the complex
 dilatation $\mu _{F_e}$ is minimal.
More precisely, if $F:\D\to \D$ is quasiconformal map such that
$F|_{\p \D}=f_4$, then its Beltrami coefficient satisfies
$||\mu_F||_{L^\infty}\geq ||\mu_{F_e}||_{L^\infty}$ and 
the equality holds only if $F=F_e$.
Furthermore, the extremal $F_e$} is a Teichm\"uller mapping, i.e.,
its  complex
 dilatation $\mu _{F_e}$ is of the form
\beq\label{phi formula}
\mu _{F_e} (z) = ||\mu _{F_e}||_{\infty}\,
 \frac{\overline{\phi(z)}}{|\phi(z)|},
\eeq
where $\phi:\D\to \C$ is holomorphic in $\D$, and has thus
discrete set of zeros.
Note that as $F_e$ need not to be even Lipschitz smooth near zeros of $\phi$.
Reader should also note that certain assumption 
on the regularity of the boundary 
value $f_m$ is necessary for existence of extremal maps. This will
be discussed after finishing the proof.

Let us now consider how a quasiconformal map $F:\D\to \D$ with
complex dilatation  $\mu_F$
 change maximal
anisotropy of conductivities. When $\sigma$ is an
isotropic conductivity in $\D$, that is, $K(\sigma)=0$,
one sees that for the anisotropic conductivity
 $\tilde \sigma=F_*\sigma$ we have
\ba
K(x,\tilde \sigma)=\mu_F(F^{-1}(x)),\quad \hbox{for $x\in \overline\D,$}
\ea 
and hence the maximal anisotropy satisfies
 $K(\tilde \sigma)=||\mu_F||_{L^\infty}$.

Let now $\gamma_3=(F_3)_*\gamma_0$ be an isotropic conductivity
in $\D$ and let  $\sigma_e=(F_e)_*\gamma_3$
be an anisotropic conductivity in $\D$.
Here, $F_e\circ F_3\circ f_m|_{\p \Omega}=f_m$.
In particular, the above {\newtext shows that
\ba
(f_m)_*\Lambda_\gamma=(f_4\circ f_3\circ f_m)_*\Lambda_{\gamma}
=(f_4\circ f_3)_*\Lambda_{\gamma_0}=(f_4)_*\Lambda_{\gamma_3}=
\Lambda_{\sigma_e}.
\ea
In particular, this implies that inverse problem
of finding conductivities $\sigma$ in $\D$
such that $(f_m)_*\Lambda_\gamma=\Lambda_{\sigma}$
has a solution $\sigma=\sigma_e$.} By Section \ref{sect: 2},
the knowledge of the boundary $\p \Omega_m=\p \D$
and the map $(f_m)_*\Lambda_\gamma$ determines
the class $\Sigma(\sigma_e)$ of conductivities in $\D$. Now
we can write the class $\Sigma(\sigma_e)$
also as
\ba
\Sigma(\sigma_e)=\{ F_*\gamma_3\ : \ 
F:\D\to \D \hbox{ is homeomorphism, }F,F^{-1}\in H^1(\Omega;\C), F|_{\p \D}=f_4\}.
\ea
Since
\ba
K(F_*\gamma_3)=||\mu_F||_{L^\infty(\D)},
\ea
we see that the conductivity $\sigma_e=(F_e)_*\gamma_3$
corresponding to the extremal map $F_e$ 
is the unique conductivity $\sigma$ in the class
$\Sigma(\sigma_e)$ that has the smallest possible value
of $K(\sigma)$.

Finally, since $|\mu _{F_e} (z)|=c_0$ is constant function 
of $z\in \D$, and $\sigma_e=(F_e)_*\gamma_3$
{\newtext with isotropic $\gamma_3$}, we see that
the ratio of the eigenvalues of the conductivity matrix
$\sigma_e(z)$ is constant for $z\in \D$.
Thus $\sigma_e$ has the form
$\sigma_e=\hat \gamma_{\lambda,\theta,\eta}$
with $c_0=(1-\lambda)/(1+\lambda)$, $\eta=\gamma_3\circ (F_e)^{-1}$,
and some $\theta$.
This proves Theorem \ref{seainoaoikea}.
\hfill $\Box$
\medskip

As noted above, the construction
of $\hat \gamma$ in the above proof
did not use the fact that $\gamma$ is isotropic  at all.
This shows Remark 1.

\medskip

Next we proof Proposition \ref{ainoatoikeat}.

{\bf Proof.} Consider isotropic conductivities $\gamma_1$ and
$\gamma_2$ in $\Omega$. 
In sequel, we use notation of the proof of Theorem \ref{seainoaoikea}.
By definition, $f_m$ determines a map $F_m$.
The construction of the map $F_1$ is based on the
Beltrami coefficient of the conductivity. Clearly,
the Beltrami coefficients for the conductivities $(F_m)_*\gamma_1$ and
$(F_m)_*\gamma_2$ coincide, and thus $F_1$ and $\Omega_1$
can be taken to be  the same  
for both $\gamma_1$ and $\gamma_2$. The maps $F_2$, $F_3$
and $F_e$ are constructed by using $\p \Omega_1$ and $F_1$,
and thus they coincide 
for $\gamma_1$ and $\gamma_2$. Since in general
$\det(F_*\gamma)(x)=\det(\gamma (F^{-1}(x)))$, this proves
 Proposition \ref{ainoatoikeat}.
\hfill $\Box$
\medskip

As noted above, certain assumption on the regularity of the boundary 
values are 
necessary, for there are counterexamples to the uniqueness, for example
 the so--called Strebel's 
chimney.
 The current state of the uniqueness question can be found in
 \cite{serb}.


We note also that if in formula (\ref{phi formula})
the function $\phi$ has zeros
in $\Omega$, then $\mu_{F_e}$
 has a singularity
of type $\overline {(z-z_0)^j}/(z-z_0)^j$, and this could affect the
 behaviour of the reconstruction algorithm
we propose in a way to explained later. However, in all the numerical 
examples we have tested these difficulties do not appear, 
probably since our deformations are relatively small.

\bigskip
\section{Electrode model}
In the numerical simulations below we have used so called
complete electrode model \cite{Somersalo92}, which is a certain
{\newtext finite dimensional approximation of Dirichlet-to-Neumann map. 
This model is chosen as it is an accurate model for the measurements
made in practice.} 
As noted before, in experimental measurements
one places the measurement electrodes on the boundary,
e.g., the skin of the patient, without knowing
exact parameterization of the boundary.
Thus this model is a paradigm of the case when 
the boundary is unknown.

To define the electrode model, let $e_j\subset \p \Omega$,
$j=1,\dots,J$
be disjoint open paths modelling the electrodes that are used
for the measurements. Let $u$ solve the equation
\beq
\label{poissoneq1}
& &\nabla\cdotp \gamma \nabla v=0\quad\hbox{in }\Omega,\\
\label{poissoneq2}
& &z_j\nu\cdotp \gamma \nabla v+v|_{e_j}=V_j,\\
\label{poissoneq3}
& &\nu\cdotp \gamma \nabla v|_{\p \Omega\setminus \cup_{j=1}^J e_j}=0,
\eeq
where $V_j$ are constants
 representing electric potentials on electrode $e_j$. 
This models the case where
electrodes $e_j$ having potentials $V_j$ are attached to the boundary,
$z_j$ is the contact impedance between electrode $e_j$ and the body surface, 
and the normal current outside the electrodes vanish. By \cite{Somersalo92},
 (\ref{poissoneq1}-\ref{poissoneq1}) has a solution $u\in H^1(\Omega)$.
The measurements in this model are the currents observed
on the electrodes, given by
\ba
I_j=\frac 1 {|e_j|}\int_{e_j}\nu\cdotp \gamma \nabla v(x)\,ds(x),\quad j=1,\dots,J.
\ea
Thus the electrode measurements are given by map
$E:\R^J\to \R^J$, $E(V_1,\dots,V_J)=(I_1,\dots,I_J)$.
We say that $E$ is the electrode measurement matrix
for $(\p \Omega,\gamma,e_1,\dots,e_J,$ $z_1,\dots,z_J)$.
Let $\Omega$ and $\tilde \Omega$ be $C^{1,\alpha}$-smooth domains.
We say that $f:\p \Omega\to \p \tilde \Omega$  
is length preserving on $\cup_{j=1}^Je_j$ if
 $||D f(\tau)||=1$ for $x\in  \cup_{j=1}^Je_j$
where $\tau$ is the unit tangent vector of $\p \Omega$.

\begin{prop}\label{Es equal}
Let $\Omega$ and $\tilde \Omega$ be $C^{1,\alpha}$-smooth domains
and $F:\overline\Omega\to \overline{\tilde \Omega}$
 be a $C^{1,\alpha}$-diffeomorphism, $e_j\subset \p \Omega$
be disjoint open sets, and $\gamma$ be a conductivity on $\Omega$.
Let $f=F|_{\p \Omega}$, $\tilde e_j=f(e_j)\subset \p \tilde \Omega$
and $\tilde \gamma=(F)_*\gamma$. Assume that $f$  
is length preserving on $\cup_{j=1}^Je_j$.
Then the electrode measurement matrices $E$
for $(\p \Omega,\gamma,e_1,\dots,e_J,z_1,\dots,z_J)$
and $\tilde E$
for $(\p \tilde \Omega,\tilde \gamma,\tilde e_1,\dots,\tilde e_J,z_1,\dots,z_J)$
coincide.
\end{prop}

{\bf Proof.}
We start with an invariant formulation of electrode
measurements $E$.
For this, let $R$ be the Robin-to-Neumann map given by
$
R f=\nu\cdotp \gamma \nabla u|_{\p \Omega}
$
{\newtext where $u$ is solution of
\beq \label{N 1}
& &\nabla\cdotp \gamma \nabla u=0\quad\hbox{in }\Omega\\
& &z\nu\cdotp \gamma \nabla v+\eta v|_{\p \Omega}=h, \nonumber
\eeq
where $z=z(x)$ is $C^\infty(\p \Omega)$ fuction such that
$z|_{e_j}=z_j$ and $\eta=\sum_{j=1}^J \chi_{e_j}(x)$,
where $\chi_{e_j}$ is the characteristic function of electrode $e_j$.} 
Note that if the boundary and the contact impedance
are known, the Robin-to-Neumann and the Dirichlet-to-Neumann
maps determine each other, that is, they represent
equivalent information.

Consider now
the bilinear form corresponding to liner maps $E:\R^J\times
\R^J\to \R$ and $R:H^{-1/2}(\p \Omega)\times H^{-1/2}(\p \Omega)\to \R$
given by
\ba
E[V,\tilde V]=\sum_{j=1}^J (EV)_j\tilde V_j |e_j| ,\quad
R[h,\tilde h]=\int_{\p \Omega} (Rh)\,\tilde h\,ds.
\ea
Let $S=\hbox{span}(\chi_{e_j}:\ j=1,\dots,J)\subset H^{-1/2}(\p \Omega)$
and define $M:V=(V_j)_{j=1}^J\mapsto \sum_{j=1}^J V_j\chi_{e_j}$
to be a map $M:\R^N\to S$. Then
\beq\label{discretization}
E[V,\tilde V]=R[MV,M\tilde V].
\eeq
Moreover, {\newtext for $h=MV$ with some $V\in \R^J$,} we have
\beq\label{form B}
R[h,h]&=&\int_{\p \Omega}
(u+z\nu\cdotp \gamma \nabla u)\nu\cdotp \gamma \nabla u\,ds
=\int_{ \Omega} \gamma \nabla u\cdotp \nabla u\,dx+
\int_{\p \Omega}z\,
|\nu\cdotp \gamma \nabla u|^2\,ds
\eeq
where $u$ solves (\ref{N 1}). Above the integral over $\Omega$
is invariant in coordinate deformations. Note that in the above 
formula the integral over the boundary
is not coordinate invariant.

Let $\tilde E$ be  electrode measurement matrix for $\tilde \gamma$
in $\tilde \Omega$ with electrodes $\tilde e_j=f(e_j)$
and let $\tilde R$ be the Robin-to-Neumann map
for $\tilde \gamma$ defined analogously to (\ref {N 1}). 
Since $f$ is length preserving on
the electrodes,
we see using (\ref{pushforward}) that $\nu\cdotp\gamma\nabla u(x)=\nu\cdotp\tilde\gamma\nabla \tilde u(f(x))$, for $\tilde u=u\circ F^{-1}$ and 
$x\in \p \Omega$, and thus we can see from (\ref{form B}) that 
\ba
R[h,\tilde h]=\tilde R[h\circ f^{-1},\tilde h\circ f^{-1}],
\ea
for $h,\tilde h\in H^{-1/2}(\p \Omega)$
supported in the closure of $\bigcup_{j=1}^Je_j$.
Thus for the map $\tilde M: V\mapsto\sum_{j=1}^J V_j\chi_{\tilde e_j}$ 
we have by formula (\ref{discretization}) that $\tilde E[V,\tilde V]=
\tilde R
[\tilde MV,\tilde M\tilde V].$ Combining this and 
(\ref{discretization}) we obtain
\ba
E[V,V']=\tilde E[V,V'].
\ea
In particular, this implies that the matrices
$E$ and $\tilde E$ coincide.
\hfill$\Box$
\smallskip

In particular, it the case where $\tilde \Omega$ is the model
domain $\Omega_m$ and $f=f_m:\p \Omega\to \p \Omega_m$
is the model map for the boundary, the assumption that
$f$ is length preserving on electrodes 
means 
the  very natural assumption that in electrode measurements 
the paramitrization of electrodes are known.
Then by Proposition \ref{Es equal}, the electrode model
discretization $E$ of $\Lambda_\gamma$ equals 
the corresponding discretization $\tilde E$ of $(f_m)_*\Lambda_\gamma$. 
Summarizing, the electrode measurements does not change
if we have modeled the geometry of the boundary incorrectly
but the electrodes are modeled correctly.

\bigskip

\section{Numerical examples}

The performance of the proposed method is evaluated by 
test cases with simulated EIT data.
First, in Section \ref{discsec} we briefly discuss 
the discretization and the computational methods that are used, 
and the results are then given in Section \ref{resultsec}. 

\subsection{Discretization and notation \label{discsec}}

The numerical solution of the forward model is based on the finite 
element method (FEM). 
The variational formulation and the finite element discretization 
of the electrode model (\ref{poissoneq1}-\ref{poissoneq3}) in the case
of isotropic conductivities have been previously discussed 
e.q. in \cite{kaipio00a}. 
The extension of the FEM-model to the anisotropic case 
is straightforward, 
the details
will be given in a subsequent publication.

For the functions $\eta(x)$ and $\theta (x)$ in equation (\ref{MAIN}) 
we use piecewise constant approximations that are defined on
a lattice of regular pixels. Thus, we have 
\bea
\eta = \sum_{i=1}^M \eta_i \chi_i (x),\quad
\theta = \sum_{i=1}^M \theta_i \chi_i (x) \label{gamdsk}
\eea
where $\chi_i$ is the characteristic function of the $i\ith$
pixel in the lattice.
Within the discretization (\ref{gamdsk}), 
the parameters $\eta$ and $\theta$ are identified 
with the coefficient vectors 
\bean
\eta = (\eta_1,\eta_2,\ldots,\eta_M)\tr \in \R^M \\
\theta = (\theta_1,\theta_2,\ldots,\theta_M)\tr \in \R^M 
\eean
and $\lambda$ is a scalar parameter. 
Note that as $\hat \gamma_{\lambda,\eta,\theta}=
\hat \gamma_{\lambda',\eta,\theta'}$, where $\lambda'=1/\lambda$ 
and $\theta'(x)=\theta(x)+\pi/2$, we can assume in looking the
minimizing uniformly anisotropic conductivity that $\lambda$ gets values
$\lambda>0$.

In practical EIT devices, the measurements are made such that known
currents are injected into the domain $\Omega$ through some
of the the electrodes
at $\partial \Omega$, and the corresponding 
voltages needeed to maintain these currents are measured on some of electrodes.
Often, voltages are measured only on those electrodes that are not
used to inject current. Thus, measurements made give only
partial information on the matrix $E$.
To take this in to account,
we introduce the following notation for the discretized problem. 
We assume that the EIT experiment consists of a set of $K$ 
partial voltage measurements, $\voltp{j}$, $j=1,\dots,K$.
For each measurement, consider a current pattern
$\currp{j},$ $j=1,\dots,K$ such that $\sum_{\ell=1}^J \currp{j}_\ell =0$.
Typically, the corresponding measurements
are the voltages (potential differences) between pairs of neighboring
electrodes. Let us assume that the measurement vector $\voltp{j}$ corresponding
to the current pattern $\currp{j}$ consist of $L$ voltage measurements, i.e., we
have $\voltp{j} \in \Rn{L}$.  
Thus, we write 
$\voltp{j} = P_jE^{-1}\currp{j} + \epsilon^{(j)}$,
where $E$ is the electrode measurement matrix, random vector $\epsilon^{(j)}$ models
the observation errors and $P_j:\R^{J}\to \R^{L}$ is a measurement 
operator that
maps the electrode potentials to measured voltages. 

In the inverse problem, the voltage measurements 
$\voltp{1},\voltp{2},\ldots,\voltp{K}$ 
are concatenated into a single vector
\[
V = (\voltp{1},\voltp{2},\ldots,\voltp{K})\tr \in \R^N,\quad N=KL.
\]
For the finite element based discretization of the forward problem 
$U: \R^{2M+1} \mapsto \R^N$, we use the notation
\[
U(\eta,\theta,\lambda) = 
(\forwp{1}(\eta,\theta,\lambda),\forwp{2}(\eta,\theta,\lambda),\ldots,
\forwp{K}(\eta,\theta,\lambda))\tr \in \R^N,
\]
respectively. 
{\newtext Here, $\forwp{j}(\eta,\theta,\lambda) = P_jE^{-1}(\eta,\theta,\lambda)\currp{j} \in \Rn{L}$
corresponds to partial voltage measurement with current pattern $\currp{j}$
and conductivity $\hat \gamma_{\eta,\theta,\lambda}$.}

Using the above notations, we {\newtext write the discretized and
 regularized version of our inverse problem} as 
finding minimizer of the functional
\be \label{ip1}
F(\eta,\theta,\lambda) = \| V - U(\eta,\theta,\lambda) \|^2 
+ W_\eta (\eta) + W_\theta (\theta) + W_\lambda (\lambda) , 
\quad \eta>0, \lambda >0, 
\ee
where the regularizing penalty functionals are of the form
\be
W_\eta (\eta) = \alpha_0 \sum_{i=1}^M \eta_i^2 +
\alpha_1 \sum_{i=1}^M \sum_{j\in \nbr_i}  | \eta_i - \eta_j |^2, 
\label{pen1}
\ee
\be
W_\theta (\theta) = \beta_0 \sum_{i=1}^M \theta_i^2 +
\beta_1 \sum_{i=1}^M \sum_{j\in \nbr_i}  | e^{i\theta_i} - e^{i\theta_j}|^2, 
\label{pen2}
\ee
\be
W_\lambda (\lambda) = \beta_2 \left( \log(\lambda) + \nu^{-2} \log(\lambda)^2 \right)
\label{pen3}
\ee
and $\nbr_i$ denotes the usual 4-point nearest neighborhood system 
for pixel $i$ in the lattice. 

Our objective is to minimize the functional (\ref{ip1}) by gradient based
optimization methods. Here we face the difficulty due to the positivity 
constraints. 
To take the 
positivity constraint into account we employ 
the interior point search method \cite{fiacco90}. In the interior
point search the original constrained problem (\ref{ip1})
is replaced by a sequence of augmented unconstrained problems
of the form
\be \label{ip2}
\tilde F_j (\eta,\theta,\lambda) = F(\eta,\theta,\lambda) + 
W_+^{(j)} (\eta)
\ee
where $W_{+}^{(j)} (\eta)$ is a penalty functional of the form
\be
W_{+}^{(j)} (\eta) = \xi_j \sum_{i=1}^M \frac{1}{\eta_i}
\label{intpen}
\ee
and $\{\xi_j\}$ is a sequence of decreasing positive parameters
such that $\xi_j \rightarrow 0$ as $j \rightarrow \infty$. Using 
a suitably chosen sequence of penalty functionals $W_{+}^{(j)}$, the
solutions of the unconstrained problems converge (asymptotically) to
the solution of the original constrained problem. 
The positivity constraint for $\lambda$ can be taken care with similar
techniques. However, it is our experience that the positivity constraint
was not needed for $\lambda$.  

For the minimization of the functionals (\ref{ip2}) we employ 
the Gauss-Newton optimization method with an explicit line search 
algorithm.

\subsection{Results \label{resultsec}}

In this section, we evaluate the performance of the proposed method with
three different test cases. The first test case is EIT data from an ellipse
domain $\domain$, in the second test case we consider an ellipse domain 
with a sharp cut and in the last test case the domain is a smooth Fourier
domain which has some resemblance with the cross section of the human body.
In all of these cases, we use the unit disk as the model domain $\compdomain$.  

In the simulations, we assume an EIT system 
with $J=16$ electrodes. In each of the test cases, the electrodes 
were located at approximately equally spaced positions at 
the exterior boundary $\partial \domain$ 
of the target domain $\domain$. The size of the electrodes 
were chosen such that the electrodes covered 
approximately 50\% of the boundary $\partial \domain$. 

The EIT measurements were simulated using the usual 
adjacent pair drive data acquisition
method. In the adjacent drive method,  
currents $+1$ and $-1$ are injected through two neighboring electrodes, 
say electrodes $e_{n}$ and $e_{n+1}$, and current through other electrodes is zero.
The voltages are measured between all $J$ pairs of neighboring electrodes.
However, three of these measurements are typically neglected
since they include either one or both of the current feeding 
electrodes $e_{n}$ or $e_{n+1}$. The rationale behind this is that
the electrode contact impedances $z_j$ are usually not known accurately. The 
possible errors in the contact impedance values cause a systematic
error between the measured voltage and the forward model for the 
measurement made
on the current feeding electrodes, and this error causes artefacts to 
the numerical reconstruction, see e.g. \cite{kolehmainen97e}. 
Thus, with the adjacent pair drive method each partial measurement consist of
$L=J-3$ voltage measurements and we have $\voltp{j} \in \Rn{J-3}$. This data
acquisition process is then repeated for all the $J$ pairs of adjacent electrodes,
leading to total of $N=J(J-3)$ voltage measurements for one EIT experiment. Thus, 
with the $J=16$ electrode system we have $V \in \Rn{208}$.


The simulated EIT measurements were computed using the isotropic EIT model and
the finite element method.
To simulate measurement noise, we added Gaussian random noise with standard 
deviation of 1\% of the maximum value of the simulated voltages to the data. 
In all of the following test cases we used value $z_\ell = 1$ for the electrode
contact impedances. These were assumed known in the inverse problem.

\begin{figure}[tb]
\centerline{\psfig{figure=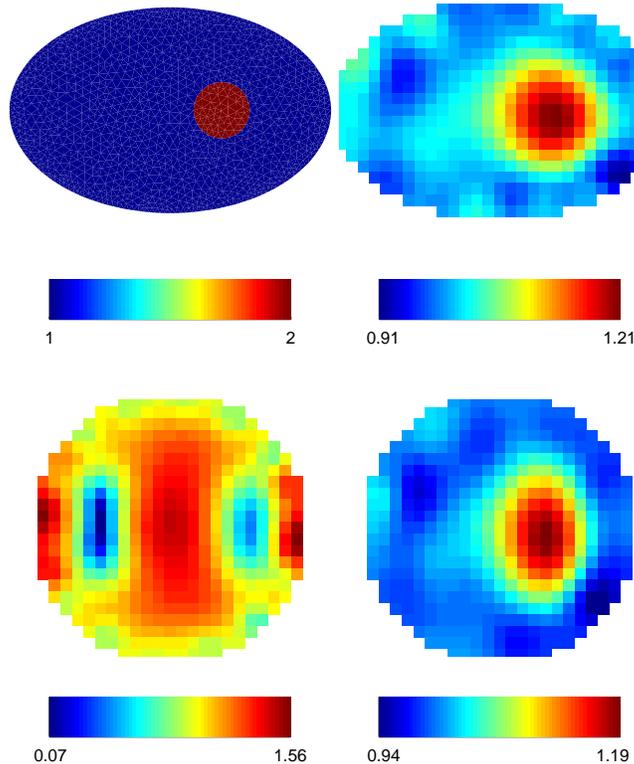,height=110mm}}
\caption{\label{res1a} Test case with EIT data from an ellipse domain
$\Omega$. The main axes of the ellipse were $1.25$ in horizontal direction 
and $0.8$ in the vertical direction. 
Top left: Simulated conductivity distribution $\gamma$. 
Top right: Reconstruction of $\gamma$ with isotropic EIT model 
in the correct domain $\domain$. Bottom left: Reconstruction
of $\gamma$ with isotropic model in incorrectly modeled 
geometry. The reconstruction domain $\compdomain$ was the unit disk. 
Bottom right: Reconstruction of $\eta$ 
with the uniformly anisotropic model in the same unit disk geometry.
}
\end{figure}
\begin{figure}[tb]
\centerline{\psfig{figure=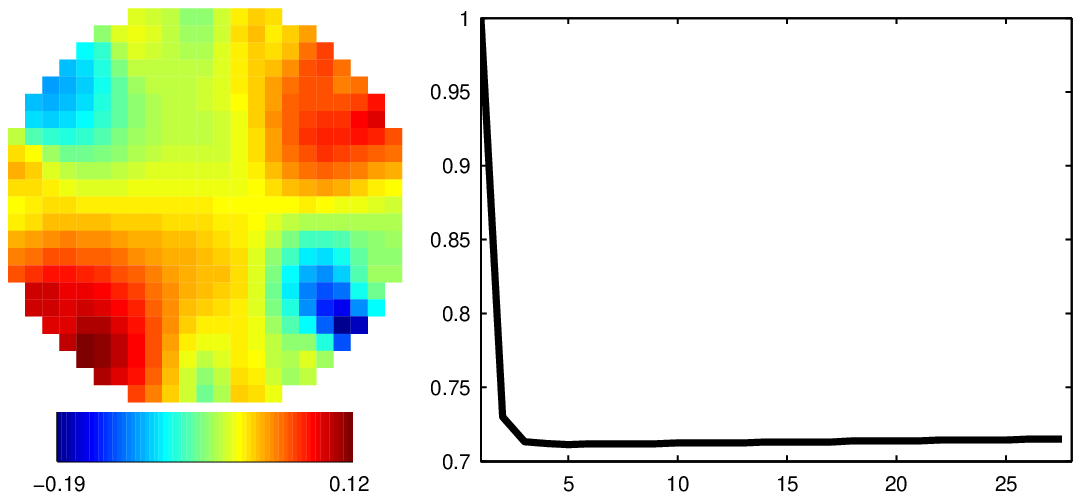,height=40mm}}
\caption{\label{res1b} Test case with EIT data from an ellipse domain
$\Omega$. The main axes of the ellipse were $1.25$ in horizontal direction 
and $0.8$ in the vertical direction. Left: Reconstruction of the 
anisotropy angle parameter $\theta$ in the incorrectly modeled geometry.
The computational domain $\compdomain$ was the unit disk. Right: Evolution
of the anisotropy parameter $\lambda$ during the Gauss-Newton iteration.}
\end{figure}

The results for the first test case are shown in Figures 
\ref{res1a}-\ref{res1b}. The target conductivity is shown in 
the top left image in Figure \ref{res1a}. The target domain 
$\domain$ is an ellipse with main axes 1.25 in the horizontal
direction and 0.8 in the vertical direction. For the simulation 
of the EIT measurements, the domain was discretized into a finite element
mesh that consisted of 1256 nodal points and 2350 triangular elements.

The reconstruction of the conductivity $\gamma$ with isotropic EIT model 
in the correct domain $\domain$ is shown
in the top right image in  Figure \ref{res1a}. The reconstruction was 
obtained by using similar optimization techniques that are explained 
{\newtext in the previous section.} 
However, in the case of isotropic model 
the unknown parameter vector is the conductivity 
vector $\gamma \in \R^M$ and the optimization 
functionals for the interior point search can be written as 
\be \label{ipiso}
H_j (\gamma) = \| V - U(\gamma) \|^2 + W_\gamma (\gamma) + W_+^{(j)} (\gamma),
\ee
where $U(\gamma)$ denotes the forward problem for the isotropic model,
$W_\gamma (\gamma)$ and $W_+^{(j)} (\gamma)$ are defined by 
equations (\ref{pen1}) and (\ref{intpen}), respectively. 
To compute the reconstruction in the top right image in  Figure \ref{res1a},
the domain $\domain$ was triangulated to a finite element mesh that consisted
of 2326 elements with 1244 nodal points. The conductivity was represented in 
a lattice of $M=451$ pixels (i.e., $\gamma \in \R^{451}$).
The regularization parameters for the 
penalty functional $W_\gamma (\gamma)$ in equation (\ref{ipiso}) 
were $\alpha_0 = 10^{-8}$ and $\alpha_1 = 10^{-4}$. When computing the 
reconstruction in the correctly modeled 
geometry, the interior point search was 
kept inactive (i.e., the sequence $\{ \xi_j \}$ of interior point search
parameters were all zeros). The conductivity vector 
was initialized to constant value of one in the optimization process. 
The Gauss-Newton optimization algorithm was iterated until convergence 
was obtained.

The image in the bottom left in Figure \ref{res1a} shows the reconstruction
of the conductivity with the isotropic model in 
incorrectly modeled geometry $\compdomain$. 
In this case, the computational domain $\compdomain$ 
was the unit disk which was triangulated to 2190 elements with 
1176 nodal points. The conductivity parameters were represented in a lattice of
$M=437$ pixels (i.e., $\gamma \in \R^{437}$). 
The regularization parameters for the penalty functional $W_\gamma (\gamma)$ 
in equation (\ref{ipiso}) 
were $\alpha_0 = 10^{-8}$ and $\alpha_1 = 2 \cdot 10^{-4}$. The sequence of 
interior point search parameters $\{\xi_j\}$ were from $2\cdot10^{-5}$ to
$5\cdot10^{-6}$. The constant vector $\gamma = 1 \in \R^{437}$ 
was used as the initial guess in the Gauss-Newton optimization.

The image in the bottom right in Figure \ref{res1a} shows the reconstruction
of $\eta$ with the uniformly anisotropic model in incorrectly modeled geometry
$\compdomain$. Here, by the solution in {\em the uniformly anisotropic model} 
we mean the optimal solution of the form (\ref{MAIN}) of the minimization
problem. 
The reconstruction was obtained by minimizing
sequence of optimization functionals of the form (\ref{ip2}). 
The reconstructed angle parameter $\theta$ is shown in the left 
image in Figure \ref{res1b}, and the evolution of the parameter 
$\lambda$ during the iteration 
is shown in the right image in Figure \ref{res1b}. 
The computational domain $\compdomain$ was the unit disk.
The finite element triangularization and the amount of the 
image pixels were the same as in the isotropic case 
in bottom left image in Figure \ref{res1a}. Thus,
the unknowns in the inverse problem are $\eta \in \R^{437}$, 
$\theta \in \R^{437}$ and $\lambda \in \R$. 
The parameters for the regularizing penalty functionals 
$W_\eta(\eta)$ in equation (\ref{pen1}) were 
$\alpha_0 = 10^{-8}$ and $\alpha_1 = 10^{-4}$. 
The parameters for the penalty functionals
$W_\theta(\theta)$ and $W_\lambda(\lambda)$ in equations 
(\ref{pen2}-\ref{pen3}) were $\beta_0 = 10^{-8}$, $\beta_1 = 5\cdot 10^{-6}$
and $\beta_2 = 0$, respectively. The sequence of 
interior point search parameters $\{\xi_j\}$ was 
from $1\cdot10^{-5}$ to $1\cdot10^{-12}$.
The Gauss-Newton optimization was started from
the constant values $\eta = 1 \in \R^{437}$, $\theta=0\in \R^{437}$ 
and $\lambda=1$ which correspond to isotropic unit conductivity.

\begin{figure}[htb]
\centerline{\psfig{figure=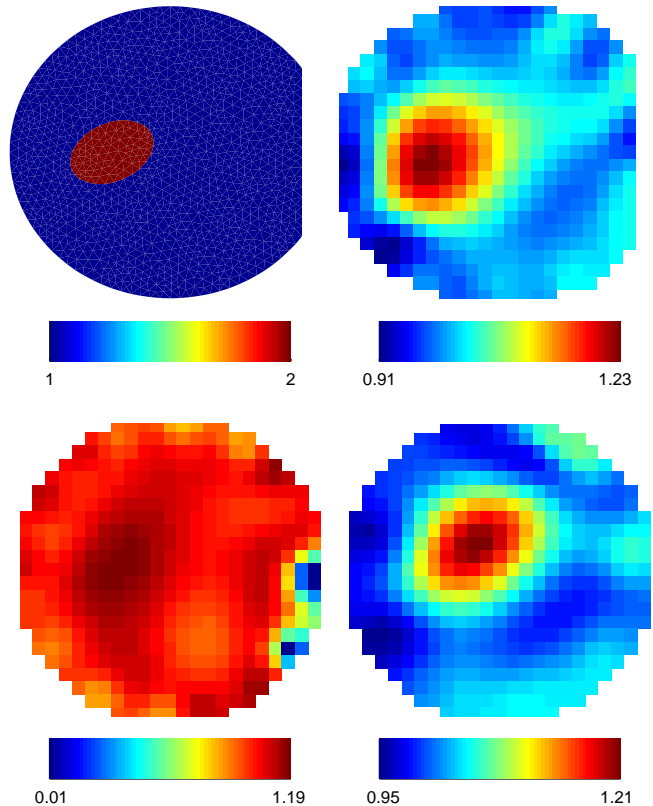,height=110mm}}
\caption{\label{res2a} Test case with EIT data 
from a truncated ellipse domain $\Omega$ with main axes $a=1.1$ and $b=0.9$.
Top left: Simulated conductivity
distribution $\gamma$. 
Top right: Reconstruction of the conductivity $\gamma$ with 
isotropic model in the correct geometry $\Omega$. 
Bottom left: Reconstruction of $\gamma$ 
with the isotropic model in incorrectly modeled geometry. 
The reconstruction domain $\compdomain$ was the unit disk. 
Bottom right: Reconstruction of the parameter $\eta$ 
with the  uniformly anisotropic model in the same unit disk geometry.}
\end{figure}

The results for the second test case are shown in Figure \ref{res2a}. The 
simulated conductivity distribution is shown in the top left image. 
In this case the domain $\domain$ is a truncated ellipse with main axes 
1.1 in the horizontal direction 
and 0.9 in the vertical direction, respectively. 
For the simulation of the EIT measurements, the 
domain was divided to a finite
element mesh of 2383 triangular elements with 1240 nodes. 

The top right image in Figure \ref{res2a} shows 
the reconstruction of the conductivity
with the isotropic model in the correct geometry. 
For the reconstruction, 
the domain $\domain$ was divided to a finite element 
mesh of 2337 triangular elements with 1217 nodes and 
the conductivity was represented on a lattice of $M=455$ pixels. 
Thus, the unknown parameter vector was $\gamma \in \R^{455}$. 
The regularization parameters for the penalty functional 
$W_\gamma (\gamma)$ in equation (\ref{ipiso}) 
were $\alpha_0 = 10^{-8}$ and $\alpha_1 = 10^{-4}$. The sequence of 
interior point search parameters $\{\xi_j\}$ was all zeros. 
The Gauss-Newton optimization was started from the constant unit conductivity.

The bottom left image in Figure \ref{res2a} shows the reconstructed
conductivity with the isotropic model in the incorrectly modeled 
geometry. 
The reconstruction
domain $\compdomain$ was the unit disk. The finite element 
mesh and pixel lattice were the same that were used for the unit disk in
Figure \ref{res1a}. Thus, the unknown conductivity vector was 
$\gamma \in \R^{437}$. The parameters for the regularizing penalty functional
$W_\gamma (\gamma)$ were $\alpha_0 = 10^{-8}$ and $\alpha_1 = 10^{-4}$, and
the sequence of interior point search parameters $\{\xi_j\}$ was from
$10^{-5}$ to $10^{-8}$. The constant unit conductivity was used as the initial
guess in the optimization. 
 
The bottom right image in Figure \ref{res2a} shows the reconstruction
of $\eta$ with the  uniformly anisotropic model in the incorrectly modeled geometry. 
The computational domain $\compdomain$ was the unit disk with the same
discretization that was used in Fig. \ref{res1a}.
Thus, the unknowns were $\eta \in \R^{437}$, $\theta \in \R^{437}$ 
and $\lambda \in \R$. 
The parameters for the regularizing penalty functionals 
$W_\eta(\eta)$ in equation (\ref{pen1}) were 
$\alpha_0 = 10^{-8}$ and $\alpha_1 = 10^{-4}$. 
The parameters for the penalty functionals
$W_\theta(\theta)$ and $W_\lambda(\lambda)$ in equations 
(\ref{pen2}-\ref{pen3}) were $\beta_0 = 10^{-8}$, $\beta_1 = 5\cdot 10^{-6}$
and $\beta_2 = 0$, respectively. The sequence of parameters 
$\{\xi_j\}$ was from $10^{-5}$ to $10^{-12}$.
The initializations for the parameters in the Gauss-Newton optimization 
were the constant values $\eta = 1 \in \R^{437}$, $\theta=0\in \R^{437}$ 
and $\lambda=1$. 

\begin{figure}[tb]
\centerline{\psfig{figure=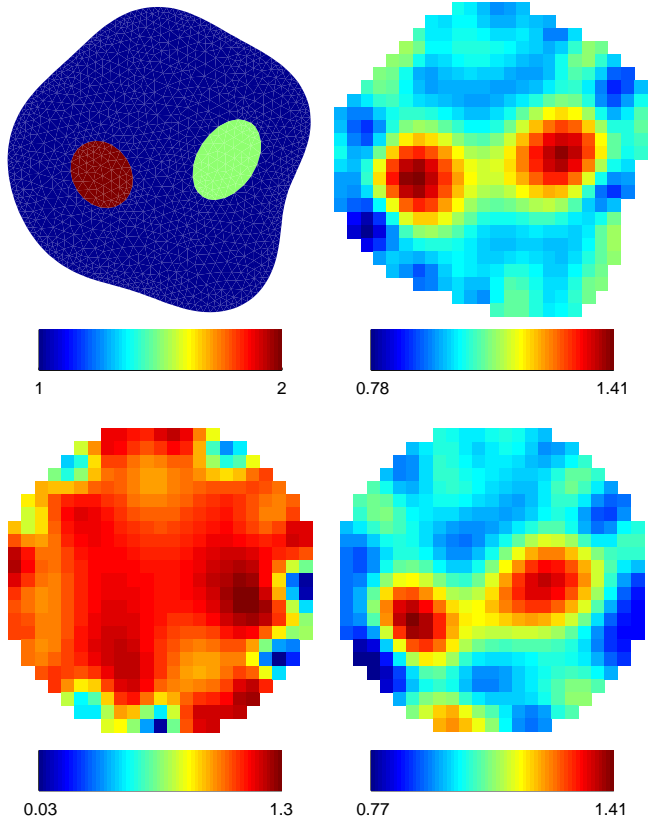,height=110mm}}
\caption{\label{res3a} 
Test case with EIT data 
from an arbitrary domain $\Omega$.
Top left: True conductivity distribution $\gamma$. 
Top right: Reconstruction of the conductivity $\gamma$ with 
isotropic model in the correct geometry $\Omega$. 
Bottom left: Reconstruction of $\gamma$ 
with the isotropic model in incorrectly modeled geometry. 
The reconstruction domain $\compdomain$ was the unit disk. 
Bottom right: Reconstruction of the parameter $\eta$ 
with the  uniformly anisotropic model in the same unit disk geometry.}
\end{figure}

The results for the last test case are shown in Figure \ref{res3a}.
In this case, the target domain $\domain$ is bounded by a smooth Fourier
boundary $\partial \domain$. 
The true isotropic conductivity distribution within the  
domain $\domain$ is shown in the top left image in Figure \ref{res3a}. 
For the simulation of the EIT measurements, the domain $\domain$ was
divided to a mesh of 2316 triangular elements with 1239 nodes.

The reconstruction of the conductivity $\gamma$ with the isotropic model
in the correct geometry $\domain$ is shown in the top right image in Figure
\ref{res3a}. The domain was divided to a mesh of 2200 triangular elements
with 1181 nodes for the image reconstruction process. The number of pixels
was $M=446$ for the representation of the conductivity image 
(i.e. $\gamma \in \R^{446}$). 
The regularization parameters for the penalty functional 
$W_\gamma (\gamma)$ were $\alpha_0 = 10^{-8}$ and $\alpha_1 = 10^{-5}$ and 
the sequence of parameters $\{\xi_j\}$ was all zeros. The constant unit 
conductivity was used as the initial guess in the Gauss-Newton optimization
algorithm.

The reconstruction of the conductivity $\gamma$ with the isotropic model
in the incorrectly modeled 
geometry is shown in the bottom left image in Figure
\ref{res3a}. The reconstruction domain $\compdomain$ was the unit disk.
The finite element mesh and the pixel lattice were the same that were used
in Figures \ref{res1a}-\ref{res2a}. Thus, the parameter vector in the inverse
problem was $\gamma \in \R^{437}$.   
The parameters in the penalty functional $W_\gamma (\gamma)$ 
were $\alpha_0 = 10^{-8}$ and $\alpha_1 = 2\cdot 10^{-4}$, and
the sequence of parameters $\{\xi_j\}$ was from
$2 \cdot 10^{-5}$ to $5 \cdot 10^{-6}$. 
The constant unit conductivity was used as the initial
guess in the optimization. 

The reconstruction of $\eta$ with 
the  uniformly anisotropic model in the incorrectly modeled geometry
is shown in the bottom right image in 
Figure \ref{res3a}. The reconstruction domain $\compdomain$ 
was the unit disk with the same discretization 
as in Figures. \ref{res1a}-\ref{res2a}.
Thus, the unknown 
parameter vectors were $\eta \in \R^{437}$, $\theta \in \R^{437}$ 
and $\lambda \in \R$. 
The parameters for the regularizing penalty functionals 
$W_\eta(\eta)$ in equation (\ref{pen1}) were 
$\alpha_0 = 10^{-8}$ and $\alpha_1 = 10^{-5}$. 
The parameters for the penalty functionals
$W_\theta(\theta)$ and $W_\lambda(\lambda)$ in equations 
(\ref{pen2}-\ref{pen3}) were $\beta_0 = 10^{-8}$, $\beta_1 = 5\cdot 10^{-6}$
and $\beta_2 = 0$, respectively. The sequence of parameters 
$\{\xi_j\}$ was from $10^{-5}$ to $10^{-12}$.
The initializations for the image parameters were the constant 
values $\eta = 1 \in \R^{437}$, $\theta=0\in \R^{437}$ 
and $\lambda=1$. 

\section{Discussion}
As can be seen from Figures \ref{res1a}-\ref{res3a}, the proposed approach
gives good results. In all test cases, the traditional 
reconstructions with the isotropic model are erroneous when the imaging
geometry is modeled incorrectly. 
The effects of erroneous geometry are seen in the
reconstructions as distortions and severe artefacts, 
especially near the boundary.  
On the other hand, 
the reconstructions of $\eta$ with the  uniformly anisotropic model 
in the same erroneous geometry are clear of these artefacts 
and represent a deformed picture of the original isotropic conductivity.
These results indicate that the proposed method offers an efficient tool
to eliminate the difficulties that arise from inaccurately known geometry
in practical EIT experiments.
 

\section*{Acknowledgements}
The authors are thankful for prof.~Kari Astala and prof.~Seppo Rickman for discussions on quasiconformal maps that were crucial for obtained results. This work was supported by the Academy of Finland (projects 203985,
72434, and 102175).

\end{document}